\documentclass[a4paper,11pt,reqno]{amsart}
\usepackage{amsmath}
\usepackage{amssymb}
\usepackage{xcolor}
\usepackage{enumerate}

\usepackage{mathrsfs}

\theoremstyle{plain} 
\newtheorem{theorem}{Theorem}[section]

\newtheorem{lemma}[theorem]{Lemma}

\makeatletter
    
    \@addtoreset{equation}{section}
  \makeatother

\title[Two independent elephant random walks]{How often can two independent elephant random walks on $\mathbb{Z}$ meet?} 
\thanks{R.R. thanks Keio University for its hospitality during multiple visits. M.T. and H.T. thank the Indian Statistical Institute for its hospitality. M.T. is partially supported by JSPS KAKENHI Grant Numbers JP19H01793, JP19K03514 and JP22K03333. H.T. is partially supported by JSPS KAKENHI Grant Number JP19K03514, JP21H04432 and JP23H01077. 
}
\author{Rahul Roy}
\address{Indian Statistical Institute, New Delhi, India}
\email{rahul@isid.ac.in}
\author{Masato Takei}
\address{Department of Applied Mathematics, Faculty of Engineering, Yokohama National University, Yokohama, Japan}
\email{takei-masato-fx@ynu.ac.jp}
\author{Hideki Tanemura}
\address{Department of Mathematics, Keio University, Yokohama, Japan}
\email{tanemura@math.keio.ac.jp}

\begin{document}

\begin{abstract}
We show that two independent elephant random walks on the integer lattice $\mathbb{Z}$ meet each other finitely often  or infinitely often depends on whether the memory parameter $p$ is strictly larger than $3/4$ or not. Asymptotic results for the distance between them are also obtained. 
\end{abstract}

\maketitle

\section{Introduction and results}

Over $100$ years have passed since the term ``random walk" appeared in the letter of Pearson \cite{Pearson1905}. It is arguably the simplest stochastic process. P\'{o}lya \cite{Polya1919RW,Polya1921RW} started the study of random walk on the $d$-dimensional lattice: A wandering point (or a random walker) moves uniformly at random between the nearest-neighbour sites on the integer lattice $\mathbb{Z}^d$ with each step being taken independent of past steps. As is explained in P\'{o}lya \cite{Polya70incidents}, his main motivation 
was to understand whether the probability that two independent random walkers eventually meet each other is one or not.
Problems on random walkers {\it with memory}, namely their future evolution depend on the entire history of the process, appeared more recently in connection with several applications. Such examples include self-avoiding random walks, reinforced random walks, and elephant random walks. See Hughes \cite{Hughes95Vol1}, R\'{e}v\'{e}sz \cite{Revesz13Book}, Laulin \cite{Laulin22PhD}, and the references therein.

A similar problem described in the previous paragraph arises also for two independent random walkers with memory.
In this note we give an answer for the elephant random walk, introduced by Sch\"{u}tz and Trimper \cite{SchutzTrimper04}.
The first step $X_1$ of the walker is $+1$ with probability $s$, and $-1$ with probability $1-s$. For each $n=1,2,\cdots$, let $U_n$ be uniformly distributed on $\{1,\cdots,n\}$, and
\begin{align*}
X_{n+1} &=  \begin{cases}
X_{U_n} &\mbox{with probability $p$}, \\
-X_{U_n} &\mbox{with probability $1-p$}, \\
\end{cases}
\end{align*} 
where $\{U_n : n=1,2,\cdots\}$ is an independent family of random variables. The sequence $\{X_i\}$ generates a one-dimensional random walk $\{S_n\}$ by
\[ S_0:=0,\quad \mbox{and} \quad S_n= \sum_{i=1}^n X_i \quad \mbox{for $n=1,2,\cdots$.} \]
In the case $p=1/2$, $\{S_n\}$ is essentially the usual symmetric random walk. For $p>1/2$ [resp. $p<1/2$] the walker prefers to do the same as [resp. the opposite of] the previous decision.

Sch\"{u}tz and Trimper \cite{SchutzTrimper04} show that there are two distinct (diffusive/super-diffusive) regimes about the asymptotic behavior of the mean square displacement. After their study, 
several limit theorems describing the influence of the memory parameter $p$ 
have been studied by many authors  \cite{BaurBertoin16,Bercu18,Collettietal17a,Collettietal17b,GuerinLaulinRaschel23,KubotaTakei19JSP,Qin23}: \\
(a) When $0<p<3/4$ the elephant random walk is diffusive, and the fluctuation is Gaussian: 
\begin{align*}
 \dfrac{S_n}{\sqrt{n}} \stackrel{\text{d}}{\to} N\left(0,\dfrac{1}{3-4p}\right) \quad \mbox{as $n \to \infty$.}
\end{align*}
where $\stackrel{\text{d}}{\to}$ denotes the convergence in distribution, and $N(\mu,\sigma^2)$ is the normal distribution with mean $\mu$ and variance $\sigma^2$. \\
(b) When $p=3/4$ the walk is marginally superdiffusive and 
\begin{align*}
 \dfrac{S_n}{\sqrt{n \log n}} \stackrel{\text{d}}{\to} N(0,1) \quad \mbox{as $n \to \infty$.}
\end{align*}
(c) If $3/4<p<1$ then there exists a random variable $L$ with a continuous distribution such that
\begin{align}
\dfrac{S_n}{n^{2p-1}} \to L \quad \mbox{a.s. and in $L^2$ as $n \to \infty$.}
\label{conv:M_nL2} \tag{$\star$}
\end{align}
Although $L$ is non-Gaussian, the fluctuation from the ``random drift" $Ln^{2p-1}$ due to the strong memory effect is still Gaussian:
\begin{align*}
\dfrac{S_n-Ln^{2p-1}}{\sqrt{n}} \stackrel{\text{d}}{\to} N\left(0,\dfrac{1}{4p-3}\right)\quad \mbox{as $n \to \infty$.}
 \label{eq:ERWCLT>1/2}
\end{align*}

Our main result is the following.

\begin{theorem} \label{thm:Main}
If $0 < p \leq 3/4$ then two independent elephant random walks with the same memory parameter $p$ meet each other infinitely often with probability one. On the other hand, if $3/4 < p <1$ then they meet each other only finitely often with probability one. 
\end{theorem} 

This theorem shows that two elephant random walks cannot meet infinitely often if the memory effect is too strong.

Theorem \ref{thm:Main} follows from

\begin{theorem} \label{thm:Quantitative}
Let $\{S_n\}$ and $\{S'_n\}$ be two independent elephant random walks with the same memory parameter $p$. 
\begin{itemize}
\item[(i)] If $0<p<3/4$ then
\begin{align*}
\limsup_{n \to \infty} \pm \dfrac{S_n - S'_n}{\sqrt{n \log \log n}} = \dfrac{2}{\sqrt{3-4p}}\quad \mbox{a.s..}
\end{align*}
\item[(ii)] If $p=3/4$ then
\begin{align*}
\limsup_{n \to \infty} \pm \dfrac{S_n - S'_n}{\sqrt{n\log n \log \log \log  n}} = 2\quad \mbox{a.s..}
\end{align*}
\item[(iii)] If $3/4<p<1$ then
\begin{align*}
\lim_{n \to \infty} \dfrac{S_n - S'_n}{n^{2p-1}} = M \quad \mbox{a.s.,}
\end{align*}
where $M$ is a random variable with
\[ P(M \neq 0)=1. \]
\end{itemize}
\end{theorem}

\section{Proof of Theorem \ref{thm:Quantitative}}

We use the following strong approximation result.

\begin{lemma}[Coletti, Gava and Sch\"{u}tz \cite{Collettietal17b}] \label{lem:StrongApprox} Let $\{S_n\}$ be the elephant random walk with the memory parameter $p$, and $\{B(t)\}$ be the standard Brownian motion.
\begin{itemize}
\item[(i)] If $0<p<3/4$ then
\begin{align*}
S_n - \dfrac{n^{2p-1}}{\sqrt{3-4p}} \cdot B(n^{3-4p}) 
&= o(\sqrt{n \log \log n}) \quad \mbox{a.s..}
\end{align*}
\item[(ii)] If $p=3/4$ then
\begin{align*}
S_n - \sqrt{n} \cdot B(\log n) 
&= o(\sqrt{n \log n \log \log \log n}) \quad \mbox{a.s..}
\end{align*}
\end{itemize}
\end{lemma}

Let $\{B(t)\}$ and $\{B'(t)\}$ be two independent standard Brownian motions. It is straightforward to see that $\{B(t)-B'(t)\}$ and $\{\sqrt{2} B(t)\}$ have the same distribution. By the law of the iterated logarithm for the standard Brownian motion,
\begin{align*}
\limsup_{t \to \infty} \dfrac{B(t) - B'(t)}{\sqrt{t \log \log t}} = 2
\end{align*}
with probability one. Thus we have 
\begin{align*}
\limsup_{n \to \infty} \dfrac{n^{2p-1} \{ B(n^{3-4p}) - B'(n^{3-4p}) \}}{\sqrt{n \log \log n}} = 2 
\end{align*}
and
\begin{align*}
\limsup_{n \to \infty} \dfrac{\sqrt{n} \{ B(\log n) - B'(\log n) \}}{\sqrt{n \log n \log \log \log n}} = 2 
\end{align*}
with probability one. Now Theorem \ref{thm:Quantitative} (i) [resp. (ii)] follows from  Lemma \ref{lem:StrongApprox} (i) [resp. (ii)].

Now we turn to the case $3/4<p<1$. By \eqref{conv:M_nL2}, 
\begin{align*}
\lim_{n \to \infty} \dfrac{S_n}{n^{2p-1}} = L \quad \mbox{and} \quad 
\lim_{n \to \infty} \dfrac{S'_n}{n^{2p-1}} = L' \quad \mbox{a.s..} 
\end{align*}
Noting that $L$ and $L'$ are independent, and both of them have continuous distributions, we have that $P(L = L')=0$. Putting $M:=L-L'$, we obtain Theorem \ref{thm:Quantitative} (iii).
\qed

\end{document}